\documentclass[11pt]{article}

\RequirePackage[OT1]{fontenc}

\RequirePackage{imsart}
\usepackage{amssymb,amsfonts,latexsym,chicago, epsf, here, graphicx,pictexwd}

\startlocaldefs

\newtheorem{lemma}{Lemma}[section]

\newtheorem{prop}{Proposition}[section]

\newtheorem{defn}{Definition}[section]

\newtheorem{remark}{Remark}[section]

\newtheorem{conj}{Conjecture}[section]

\newtheorem{example}{Example}[section]

\newcommand{\RR}{ {\mathbb R}}

\endlocaldefs

\begin{document}

\begin{frontmatter}

\title{Conjecture of error boundedness in a \\
new Hermite interpolation problem via splines of odd-degree}
\runtitle{Conjecture:  error boundedness Hermite interpolation}


\author{\fnms{Fadoua} \snm{Balabdaoui}\thanksref{t1,t3}\corref{}\ead[label=e1]{fadoua@math.uni-goettingen.de}}
\thankstext{t1}{Research supported in part by National Science Foundation
grant DMS-0203320}
\address{Institute for Mathematical Stochastics\\
Georgia Augusta University Goettingen\\
Maschmuehlenweg 8-10\\
D-37073 Goettingen\\
Germany \\\printead{e1}}
\and
\author{\fnms{Jon A} \snm{Wellner}\thanksref{t2}\ead[label=e2]{jaw@stat.washington.edu}}
\thankstext{t2}{Research supported in part by National Science Foundation grant
DMS-0203320,  NIAID grant 2R01 AI291968-04, and an NWO Grant
to the Vrije Universiteit, Amsterdam}
\thankstext{t3}{Corresponding author}
\address{Department of Statistics\\Box 354322 \\University of Washington\\Seattle, WA 98195-4322 \\
\printead{e2}}

\affiliation{University of G\"ottingen and University of Washington}

\runauthor{Balabdaoui and Wellner}

\begin{abstract}
We present a Hermite interpolation
problem via splines of odd-degree which, to the best knowledge
of the authors, has not been considered in the literature on interpolation
via odd-degree splines. In this new interpolation problem, we conjecture
that the interpolation error is bounded in the supremum norm
\textit{independently} of the locations of the knots.
Given an integer $k \geq 3$, our spline interpolant is of degree $2k-1$ and
with $2k-4$ (interior) knots.
Simulations were performed to
check the validity of the conjecture.
We present strong
numerical evidence in support of the conjecture for $k=3, \cdots, 10$ when the interpolated
function belongs to $C^{(2k)}[0,1]$, the class of $2k$-times continuously
differentiable functions on $[0,1]$. In this case, the worst interpolation
error is proved to be attained by the perfect spline of degree $2k$
with the same knots as the spline interpolant.
This interpolation problem arises naturally in nonparametric
estimation of a multiply monotone density via Least Squares and Maximum Likehood methods.
\end{abstract}



\end{frontmatter}

\section{Introduction}

Interpolation via splines of degree $2k-1$ with simple knots has been
considered for different boundary conditions that the spline interpolant
is required to satisfy (see e.g. [Schoenberg 1964a; 1964b] and
[N\"urnberger 1989, pages 116-123]. One particular aspect of these
problems is that the $k$-th derivative of the spline interpolant gives
the minimal $L_2$-norm among all functions that are smooth enough
and satisfy the same boundary conditions (see e.g. [de Boor 1964; 1974],
[Holladay 1957], [Schoenberg 1964a; 1964b; 1964c], [N\"urnberger 1989]).
In the particular case of \textit{complete} and \textit{natural} interpolation
of elements in the Sobolev space $W^l_p$ via a spline of a given degree
$2k-1$, this optimality property was a good starting point to prove that
the the $L_p$-norm of the interpolation error is independent of the locations
of the knots (see [Shadrin 1992] and the references given there).
We recall here that given a function $f$ on $[0,1]$,
$0= \tau_0 < \tau_1 < \ldots < \tau_{m-1} < \tau_m= 1$ and $k \geq 1$ an integer,
a spline $s$ with degree $2k-1$ and (simple) knots $\tau_1, \cdots, \tau_{m-1}$
is the complete spline if and only if $\Vert s^{(k)} \Vert_2$ minimizes
$\Vert g^{(k)} \Vert_2$ over all functions $g \in W^k_2$ satisfying
\begin{eqnarray}\label{ComplCond}
&& g(\tau_j) = f(\tau_j), \ \ \textrm{for $j=0, \ldots, m$} \nonumber \\
&& g^{(i)}(\tau_0) = g^{(i)}(\tau_0), \ g^{(i)}(\tau_{m})
= f^{(i)}(\tau_m) \ \ \textrm{for all $i=1, \ldots, k-1$},
\end{eqnarray}
and is the natural spline if and only if $\Vert s^{(k)} \Vert_2$
minimizes $\Vert g^{(k)} \Vert_2$ over all functions $g \in W^k_2$ satisfying
\begin{eqnarray}\label{NaturCond}
g(\tau_j) = f(\tau_j), \ \ \textrm{for $j=0, \ldots, m$},
\end{eqnarray}
(see e.g. [Holladay 1957], [de Boor 1964; 1974],
[Schoenberg 1964a; 1964b; 1964c], [Shadrin 1992]).

Natural cubic splines ($k=2$) go back at least to [Holladay 1957].
[Holladay 1957] was apparently the first to prove that the unique
spline that satisfies the conditions in (\ref{NaturCond}), and such
that it is a cubic polynomial to
the left and to the right of $\tau_0$ and $\tau_m$ respectively is
equal to the natural spline interpolant. In general, a natural
spline $s$ of degree $2k-1 \geq 3$ admits a polynomial
extension of degree $2k-1$ to the left and right of the boundary
points $\tau_0$ and $\tau_m$ respectively; i.e.,
\begin{eqnarray}\label{DerivNatur}
s^{(i)}(\tau_0) = s^{(i)}(\tau_m)=0, \ \ \textrm{for all $i=k, \ldots, 2k-2$}
\end{eqnarray}
(see e.g. [Schoenberg 1964a; 1964b; 1964c], [N\"urnberger 1989]).

\medskip

If we fix the number of knots, and take $h$ to be the maximal
distance between two successive knots, it follows from Theorem 6.4
of [Shadrin 1992] that there exists $c_1 > 0$ (depending only on $k$) such that
\begin{eqnarray}\label{ShadIneq}
\Vert f^{(i)} - s^{(i)} \Vert_\infty \leq c_1 \ h^{l-i} \ \Vert f^{(l)} \Vert_\infty,
\end{eqnarray}
where $f \in C^{(l)} [0,1]$,
$0 \leq i < k \leq l \leq 2k $, and $s$ is either the (corresponding)
natural or complete spline interpolant ([Shadrin 1992]
makes this assertion for $f $ in the Sobolev space
$W_{\infty}^l$, but via a standard Sobolev embedding theorem (see e.g.
[Adams and Fournier, 2003], pages 80, 85)
this implies the stated inequality).
The result provides bounds
that are independent of the locations of the knots
$\tau_1, \ldots, \tau_{m-1}$. In particular, if $f$ is
$C^{(2k)}[0,1]$, and $s$ is its complete spline interpolant then
\begin{eqnarray*}
\Vert f - s \Vert_\infty \leq c_1 h^{2k} \Vert f^{(2k)} \Vert_\infty,
\end{eqnarray*}
a bound that [de Boor 1974] had conjectured for $k > 4$.
In connection with the previous bound, [de Boor 1974]
found that it was sufficient to show boundedness of the
supremum norm of the $L_2$-projector of $C^{(k)}[0,1]$
on the space of splines of degree $k-1$ with knots
$\tau_1, \ldots, \tau_{n-1}$. In other words, de Boor
conjectured that there exists some constant $d_k > 0$
(depending only on $k$)
such that for any nontrivial function $f$ in $C^{(k)}[0,1]$
\begin{eqnarray}\label{LinfBound}
\sup_{0< \tau_1 < \cdots < \tau_{m-1} < 1} \frac{\Vert Pf \Vert_\infty}{\Vert f \Vert_\infty} \leq d_k,
\end{eqnarray}
where $Pf$ is the orthogonal projection of $f$ on the space of
splines of degree $k-1$ with knots $\tau_1, \cdots, \tau_{m-1}$,
equipped with the ordinary inner product
\begin{eqnarray*}
\langle s_1, s_2 \rangle = \int_0^1 s_1(t) s_2(t) dt,
\end{eqnarray*}
(see also [Shadrin 2001]). Apparently, de Boor had conjectured
such a result for the first time in 1972 ([de Boor 1973]). However,
the conjecture remained unsolved for more than 25 years before
[Shadrin 2001] found a proof.
The bound in (\ref{LinfBound})
combined with the techniques developed by [de Boor 1974] for bounding
interpolation error implies that there exists $d_{j,k} > 0$ depending only
on $j$ and $k$ such that
\begin{eqnarray*}
\Vert f - s \Vert \leq d_{j,k} \ h^{k+j} \ \Vert f^{(k+j)} \Vert_\infty,
\end{eqnarray*}
for all functions $f \in C^{(k+j)}, j=0, \ldots, k$. The latter can
obviously be viewed as a special case of (\ref{ShadIneq}), with
$i=0$, $l=k+j$ and $d_{j,k} = c_1$. However, the bound in (\ref{LinfBound})
makes the inequality in (\ref{ShadIneq}) valid for $i=k$ as well, and
hence it provides a stronger result. \\

\section{A new Hermite interpolation problem}

\subsection{Description of the problem} 

Let $k \geq 2$ be a fixed integer, and consider
$\tau_0=0 < \tau_1 < \cdots < \tau_{2k-4} < \tau_{2k-3}=1$.
Given a function $f$ defined on $[0,1]$ and at least differentiable at the
$2k-2$ points $\tau_i$ for $ i=0, \cdots, 2k-3$, we denote by
$H_kf$ the unique spline of degree $2k-1$ with (simple) knots
$\tau_1, \cdots, \tau_{2k-4}$ that satisfies the $4k-4$ linear conditions
\begin{eqnarray}\label{NewHerm}
(H_kf)(\tau_i) = f(\tau_i), \ \ \ \textrm{and} \ \ \
(H_kf)'(\tau_i) = f'(\tau_i) \ \ \textrm{for $i=0,\cdots, 2k-3$}.
\end{eqnarray}
In what follows, we denote by $\mathcal{E}_k(f)$ the
associated interpolation error $f-H_k f$.
Existence and uniqueness of the solution follows easily
from the \\ Schoenberg-Whitney-Karlin-Ziegler
Theorem \\ ([Schoenberg and Whitney 1953]; Theorem 3, page 529,
[Karlin and Ziegler 1966] ; or see Theorem 3.7, page 109,
[N\"urnberger 1989]; or Theorem 9.2, page 162,
[DeVore and Lorentz 1993] ).
When $k=2$, the problem reduces to complete interpolation via a
cubic polynomial. The latter is of less concern since asymptotic theory 
for the motivating statistical problem
has already been developed for this case;
see [Groeneboom, Jongbloed, and Wellner (2001a), (2001b)].
We therefore assume hereafter that $k \geq 3$.

This problem seems to be new in the sense that it has not been
considered by the literature on interpolation via splines in general.
Unlike in complete or natural interpolation, the \lq \lq smoothing\rq\rq \
boundary conditions of the spline interpolant $H_kf$ are not being absorbed
at the boundary points $\tau_0$ and $\tau_m$
with $m=2k-3$ (here we refer to the conditions
in (\ref{ComplCond}) and (\ref{DerivNatur})) but instead, are \lq\lq
re-distributed\rq \rq \ somehow evenly among the knots. Another feature
of this problem is that the number of knots of the spline interpolant
depends on its degree, which seems to be somewhat unusual (personal
communication with Nira Dyn). On the other hand, functions that are
not twice differentiable or smoother on $[0,1]$ can be very well interpolated
as long as they are differentiable at $\tau_i$ for $ i=0, \cdots, 2k-3$.
Whether the interpolation error is bounded
for these functions remains an open question, but it is
clear that they are eliminated from the scheme of complete
or natural interpolation via splines of degree $\geq 3$.

\subsection{Motivation}
The consideration of this \lq\lq non-standard\rq\rq \ Hermite
interpolation problem was originally motivated by a nonparametric
estimation problem in a certain class of shape-constrained densities.
In fact, it arises naturally in connection with the asymptotic
distribution of the Least Squares (LS) and Maximum Likelihood (ML)
estimators of a $k$-monotone density $g$ on $(0,\infty)$; i.e.,
a density that satisfies: $(-1)^l g^{(l)}$ is nonincreasing and
convex for $l=0, \cdots, k-2$. Let $X_1, \cdots, X_n$ be $n$ independent
random variables with a common $k$-monotone density $g_0$.
Let $\mathbb{G}_n$ be the empirical distribution of the $X_i$'s; i.e.
\begin{eqnarray*}
\mathbb{G}_n(x) = \frac{1}{n} \sum_{j=1}^n 1_{[X_j \leq x]}.
\end{eqnarray*}
It has been shown [Balabdaoui and Wellner 2004a]
that both the ML and LS estimators, denoted hereafter by
$\hat{g}_n$ and $\tilde{g}_n$ respectively, are splines
of degree $k-1$ with simple (random) knots.
The number and the locations of these are uniquely determined by
the minimization problem defining the estimators (and the data via ${\mathbb G}_n$),
and they can be computed
via an iterative spline algorithm
(see [Balabdaoui and Wellner 2004b]). Now, when the true $k$-monotone
density $g_0$ is assumed to be $k$-times differentiable at some
fixed point $x_0 > 0$ such that $(-1)^k g^{(k)}_0(x_0) > 0$ and $g^{(k)}_0$
is continuous in the neighborhood of $x_0$, then the number of knots
around $x_0$ increases to infinity almost surely and hence the distance
between two successive knots decreases to 0. It turns out that proving
that the stochastic order of this distance is $n^{-1/(2k+1)}$ as $n \to \infty$
is the main key to deriving the exact rates of convergence of
$\tilde{g}^{(j)}_n(x_0)$ or $\hat{g}^{(j)}_n(x_0)$, $j=0, \cdots, k-1$,
as well as the corresponding limiting distribution.

Here we choose to focus on $\tilde{g}_n$ rather than $\hat{g}_n$
since it is easier through the simple characterization of $\tilde{g}_n$
to make the link between the new Hermite problem in (\ref{NewHerm})
and the nonparametric estimation problem. It has been shown
[Balabdaoui and Wellner 2004a] that a necessary and sufficient condition
for a $k$-monotone function $\tilde{g}_n$ to be equal to the LSE is given by
\begin{eqnarray}\label{CharLSE}
\left \{
\begin{array}{ll}
\tilde{H}_n(x) \geq \mathbb{Y}_n(x), \ \ \textrm{for all $x > 0$} \\
\tilde{H}_n(x) = \mathbb{Y}_n(x), \ \ \textrm{if $(-1)^k
\tilde{g}^{(k-1)}_n(x-) < (-1)^k \tilde{g}^{(k-1)}_n(x+)$}, \\
\end{array}
\right.
\end{eqnarray}
where
\begin{eqnarray*}
\tilde{H}_n(x) = \frac{1}{(k-1)!} \int_0^x (x-t)^{k-1} \tilde{g}_n(t)dt
\end{eqnarray*}
and
\begin{eqnarray*}
\mathbb{Y}_n(x) = \frac{1}{(k-1)!} \int_0^x (x-t)^{k-1} d\mathbb{G}_n(t).
\end{eqnarray*}
In other words, $\tilde{H}_n$ and $\mathbb{Y}_n$ are the $k$-fold
and $(k-1)$-fold integrals of $\tilde{g}_n$ and $\mathbb{G}_n$
respectively, and their corresponding curves have to touch
each other at every knot of $\tilde{g}_n$. Note also that $\tilde{H}_n$
is a spline of degree $2k-1$ with the same knots as $\tilde{g}_n$.
Now, if we assume that $\tilde{g}_n$ has at least $2k-2$ knots
$\tau_0 < \cdots < \tau_{2k-3}$ in the neighborhood of $x_0$,
an event that occurs with probability $\to 1$ as $n \to \infty$,
then it follows from the characterization in (\ref{CharLSE}) that
\begin{eqnarray*}
\tilde{H}_n(\tau_i) = \mathbb{Y}_n(\tau_i) \ \ \textrm{and} \ \ \tilde{H}'_n(\tau_i) = \mathbb{Y}'_n(\tau_i)
\end{eqnarray*}
for $i=0, \cdots, 2k-3$. Hence, $\tilde{H}_n = H_k\mathbb{Y}_n$
on $[\tau_0, \tau_{2k-3}]$, an identity that plays a major role in
studying the distance between consecutive knots.
Through this identity, the Least Square esimator, which in essence
estimates the whole density, could be studied locally in the
neighborhood of the point of interest $x_0$. To be able to obtain the
stochastic order of the distance between the knots of the estimator,
arguments from empirical processes theory ([van der Vaart and Wellner 1996])
are used in combination with the current conjecture on the interpolation error.
For more technical details,
we refer to [Balabdaoui and Wellner 2004d, page 12, Lemma 2.3].
Another version of this same problem arises in showing the existence
of a solution to the natural Gaussian version of the problem 
[see Balabdaoui and Wellner 2004c].

\subsection{The conjecture}

Proving that the distance between successive knots is indeed of
the stochastic order $n^{-1/(2k+1)}$ as $n \to \infty$ still
depends on bounding $\Vert \mathcal{E}_k(f) \Vert_\infty $
independently of the locations of the knots, for any function $(k-1)$-differentiable
function $f$ such that $f^{(k-1)}$ admits a finite total variation.
For more details, see [Balabdaoui and Wellner 2004d, page 16,
Lemma 2.4]. We formulate our conjecture as follows:

\bigskip

\par \noindent
\begin{conj}
Let $\tau_0=0 < \tau_1 < \cdots < \tau_{2k-4} < \tau_{2k-3}=1$
and $u \in (0,1)$. If $f_u(t)= (t- u)^{k-1}_{+}/(k-1)!$,
then there exists $c_k > 0$ 
such that
\begin{eqnarray}\label{Conjec}
\sup_{u \in (0,1)} \sup_{0 <\tau_1 < \cdots < \tau_{2k-4} < 1}
\Vert \mathcal{E}_k(f_u) \Vert_\infty \leq c_k.
\end{eqnarray}
\end{conj}
\bigskip

If the bound in (\ref{Conjec}) is true, then we can establish the following lemma:
\bigskip

\par\noindent
\begin{lemma}
If (8) holds, then
for any function $f$ in $C^{(k+j)}[0,1],
j=0, \cdots, k$, we have
\begin{eqnarray}\label{Implic.}
\Vert \mathcal{E}_k(f) \Vert_\infty \leq \frac{c_k}{(j+1)!} \ \Vert f^{(k+j)}\Vert_\infty.
\end{eqnarray}
\end{lemma}

\medskip

\par \noindent
\textbf{Proof.}
Using Taylor expansion we
can write for all $t \in (0,1)$
\begin{eqnarray*}
f(t)
& = & f(0) + f'(0)t + \cdots + \frac{f^{(k+j-1)}(0)}{(k+j-1)!}t^{k+j-1} \\
& & \qquad + \  \frac{1}{(k+j-1)!} \int_0^1 (t - u)^{k+j-1}_+ f^{(k+j)}(u) du \\
& = & P_{k,j}(t) + \frac{1}{(k+j-1)!} \int_0^1 (t - u)^{k+j-1}_+ f^{(k+j)}(u) du .\\
\end{eqnarray*}
Since $P_{k,j}$ is a polynomial of degree at most
$2k-1$ for $j=0, \cdots, k$, it follows that
\begin{eqnarray*}
\mathcal{E}_k(f)
= \frac{1}{(k+j-1)!} \int_0^1 \left [H_k(\cdot-u)^{k+j-1}_+ - (\cdot-u)^{k+j-1}_+ \right]f^{(k+j)}(u) du,
\end{eqnarray*}
on $[0,1]$, and hence, if $j=0$, we have
\begin{eqnarray*}
\Vert \mathcal{E}_k(f) \Vert_\infty \leq c_k \Vert f^{(k)}\Vert_\infty
\end{eqnarray*}
by the conjecture. Now, for $j=1, \cdots, k$, we can apply
again Taylor expansion to the function $t \mapsto (t-u)^{k+j-1}_{+}$ ($u$ fixed) to obtain
\begin{eqnarray*}
(t-u)^{k+j-1}_{+} & = & \frac{(k+j-1)\cdots j}{(k-1)!} \int_0^1 (t - v)^{k-1}_+ (v-u)^{j-1}_{+} dv.
\end{eqnarray*}
This implies that
\begin{eqnarray*}
\lefteqn{ \bigg \Vert H_k(\cdot -u)^{k+j-1}_+
- (\cdot-u)^{k+j-1}_+ \bigg \Vert_\infty} \\
&=& \bigg \Vert \frac{(k+j-1)\cdots j}{(k-1)!} \int_0^1 \left [H_k(\cdot-v)^{k-1}_+
- (\cdot-v)^{k-1}_+ \right](v-u)^{j-1}_{+} dv \bigg \Vert_\infty \\
& \leq & (k+j-1)\cdots j \ c_k \ \int_u^1 (v-u)^{j-1}dv = \frac{(k+j-1)\cdots j}{j}c_k (1-u)^j \\
& = & \frac{(k+j-1)!}{j!} c_k (1-u)^j \, .
\end{eqnarray*}
and therefore
\begin{eqnarray*}
\Vert \mathcal{E}_k(f) \Vert_\infty \leq \frac{c_k}{(j+1)!} \ \Vert f^{(k+j)}\Vert_\infty.
\end{eqnarray*}
\hfill $\blacksquare$

\medskip
\bigskip

\par \noindent
\textbf{Remark.}
To conform more closely with the error bound given in the
natural and complete interpolation problem, we note that the
bound in (\ref{Implic.}) can be replaced by $d_{k,j} h^{k+j}$,
where $h$ denotes again the largest distance between two
successive knots, and $d_{k,j} \leq c_{k,j} 2^{(2k-3)(k+j-1)} (2k-3)$,
which follows from iterative application of the $c_r$-inequality.

\section{Computations}
\subsection{Interpolating the hinge function $f_u$}
In the absence of a theoretical proof of (8), a natural approach is to replicate the interpolation problem a number of times. The knots of the spline interpolant as well as the single knot, $u$, of the hinge function $f_u$ (the function being interpolated) are drawn independently from a uniform distribution on $[0,1]$. The procedure can be simply described as follows:

\begin{enumerate}
\item
Generate $2k-3$ independent uniform random variables on $(0,1)$, \\ $U_1, \cdots, U_{2k-3}$.
\item
Set $ (\tau_1, \cdots, \tau_{2k-4}) = (U_{(1)}, \cdots, U_{(2k-4)})$ and
$u= U_{2k-3}$, where \\ $U_{(1)}, \cdots, U_{(2k-4)}$
are the order statistics of $U_1, \cdots, U_{2k-4}$.
\item
Find the spline $H_kf_u$, and calculate the supremum
norm error \\ 
$\Vert \mathcal{E}_k(f_u) \Vert_\infty$.
\item
Store $\Vert \mathcal{E}_k(f_u) \Vert_\infty$ and repeat the previous steps.
\end{enumerate}
\bigskip

The Mathematica programs \lq EB-SinglePrint-hinge-post\rq \
and \lq EB-MC-hinge-post\rq \ used to produce the data are
provided in the supporting material; see
\smallskip

\par\noindent
\begin{quote}
www.stat.washington.edu/jaw/RESEARCH/
SOFTWARE/software.list.html
\end{quote}
\smallskip

\par\noindent
for all the programs used to produce the tables and figures.
The first program produces
the maximal interpolation error for a single run as well as the plot
of the spline interpolant and that of the corresponding error.
The Monte Carlo procedure described above can be run using
the second program which produces a statistical summary as well
as plots of the empirical distribution of the interpolation error.
Choosing a large number of digits of the random numbers, $d$,
seems to resolve some of the numerical issues that we have
encountered: When some of the knots are extremely close to
each other, the linear system becomes
very close to being singular in which case Mathematica suppresses
the corresponding outcome. This suppression due to severe ill-conditioning
of the matrices does not happen if $d$ is taken to be appropriately large.
However, we believe that other numerical instabilities do occur and
become increasingly severe as we increase the value of $k$.
For $k=3, \cdots, 10$, we performed $2 \times 10^4$ Monte Carlo
replications. The corresponding statistical summaries can be found in
Table \ref{ErrorStat}. Of course, one is interested in the maximal
interpolation error over all configurations, but it is also of some
interest to get an idea about its distribution. In these simulations
and further ones, the obtained maximal error over all knot configurations
seems to be stable and reasonably small for $k=3,4,5$. However, the
order of this maximal error seems to increase with $k$ to reach
extremely large values ($10^{13}$ for $k=10 \ !$). Note that the
values of the median and the 95\%-quantile decay steadily. In the
absence of theoretical argument that proves or disproves the conjecture,
large maximal errors would cast doubt on the validity of the conjecture
for large $k$. However, we have an evidence of numerical instabilities
that will be presented and discussed in the next subsection. There,
the interpolated function is much smoother than the hinge function
since it is assumed to be $2k$-times continuously differentiable,
nevertheless extremely large interpolation errors do also occur.
These large maximal errors can be compared with much smaller
bounds that could be obtained via the perfect splines theory.
So in this particular case, both theory and simulations are
combined to seek a stronger support for uniform boundedness of the interpolation error.
\bigskip

\begin{center}
\begin{table}[12pt]
\caption{Statistical summary of the $L_\infty$-norm of the interpolation error,
$\mathcal{E}_k(f_u)$, for $k=3,\cdots,10$ based on $10^4 + 10^4$ independent simulations. Table produced using the Mathematica program \lq EB-MC-hinge-post\rq.}
\medskip
\begin{tabular}{|c|cccccc|}
\hline
$k$ & Mean & Median & Std. dev. & $95\%$-q & $99\%$-q & max \\ \hline\hline
3 & $4.54 \times 10^{-3}$ & $2.69 \times 10^{-3}$ & $8.8 \times 10^{-3}$ & $1.47 \times 10^{-2}$ & $3.44 \times 10^{-2}$ & $3.2 \times 10^{-1}$ \\
3 & $4.55 \times 10^{-3}$ & $2.59 \times 10^{-3}$ & $1.07 \times 10^{-2}$ & $1.38 \times 10^{-2}$ & $3.43 \times 10^{-2}$ & $5.4 \times 10^{-1}$ \\
\hline
4 & $5.6 \times 10^{-4}$ & $7.92 \times 10^{-5}$ & $3.72 \times 10^{-3}$ & $1.85 \times 10^{-3}$ & $8.58 \times 10^{-3}$ & $1.8 \times 10^{-1}$ \\
4 & $7.16 \times 10^{-4}$ & $7.96 \times 10^{-5}$ & $6.27 \times 10^{-3}$ & $2 \times 10^{-3}$ & $9.13 \times 10^{-3}$ & $3.2 \times 10^{-1}$ \\
\hline
5 & $4.07 \times 10^{-4}$ & $2.68 \times 10^{-6}$ & $9.11 \times 10^{-3}$ & $3.87 \times 10^{-4}$ & $5 \times 10^{-3}$ & $6 \times 10^{-1}$ \\
5 & $2.97 \times 10^{-4}$ & $2.75 \times 10^{-6}$ & $5 \times 10^{-3}$ & $4.05 \times 10^{-4}$ & $3.59 \times 10^{-3}$ & $3.35 \times 10^{-1}$ \\
\hline
6 & $1.51 \times 10^{-2}$ & $1.20 \times 10^{-7}$ & 1.42 & $1.41 \times 10^{-4}$ & $3 \times 10^{-3}$ & $1.42 \times 10^{2}$ \\
6 & $1.42 \times 10^{-2}$ & $1.17 \times 10^{-7}$ & $9 \times 10^{-2}$ & $1.40 \times 10^{-4}$ & $2.44 \times 10^{-3}$ & $8.87$ \\
\hline
7 & $9.87 \times 10^{-2}$ & $7.73 \times 10^{-9}$ & 9.6 & $6.53 \times 10^{-5}$ & $2.38 \times 10^{-3}$ & $9.6 \times 10^2$ \\
7 & $2.06 \times 10^{2}$ & $8.9 \times 10^{-9}$ & $2.05 \times 10^5 $ & $6.75 \times 10^{-5}$ & $2.89 \times 10^{-3}$ & $2.05 \times 10^6$ \\
\hline
8 & 9.11 & $6.69 \times 10^{-10}$ & $8.03 \times 10^{2}$ & $4.36 \times 10^{-5}$ & $4.41 \times 10^{-3}$ & $7.97 \times 10^{5}$ \\
8 & $5.5 \times 10^{-1}$ & $6.62 \times 10^{-10}$ & $3.92 \times 10^{1}$ & $4 \times 10^{-5}$ & $2.83 \times 10^{-3}$ & $3.41 \times 10^{4}$ \\
\hline
9 & $4.57 \times 10^{2}$ & $7.58 \times 10^{-11}$ & $3.38 \times 10^{5}$ & $2.47 \times 10^{-5}$ & $1.12 \times 10^{-2}$ & $3.31 \times 10^{6}$ \\
9 & $1.83 \times 10^{1}$ & $6.79 \times 10^{-11}$ & $1.71 \times 10^{4}$ & $2.63 \times 10^{-5}$ & $7.36 \times 10^{-3}$ & $1.71 \times 10^{6}$ \\
\hline
10 & $1.4 \times 10^{7}$ & $1.06 \times 10^{-11}$ & $1.37 \times 10^{9}$ & $2.94 \times 10^{-5}$ & $3.82 \times 10^{-2}$ & $1.37 \times 10^{11}$ \\
10 & $1.70 \times 10^{9}$ & $1.16 \times 10^{-11}$ & $1.70 \times 10^{11}$ & $3.9 \times 10^{-5}$ & $8.95 \times 10^{-2}$ & $1.7 \times 10^{13}$ \\ \hline
\end{tabular}
\label{ErrorStat}
\end{table}
\end{center}

\subsection{Interpolation of $f \in C^{(2k)}[0,1]$}

We have proved in subsection 2.3 that provided that the conjecture
is true, the interpolation error is uniformly bounded for any function
in $C^{(k+j)}[0,1], j=0, \cdots, k$ with respect to $\Vert f^{(k+j)} \Vert_\infty$.
When $j=k$, the interpolated function $f$ is $2k$-times continuously
differentiable on $[0,1]$.
It turns out in this particular case that one can take another
route to seek evidence for uniform boundedness of the interpolation
error with respect to $\Vert f^{(2k)} \Vert_\infty$.

Indeed, the interpolation error can be bounded more explicitly by the error for interpolation of a perfect
spline, as was pointed out to us by A. Shadrin. In this particular case, for
a fixed set of knots $0 < \tau_1 < \cdots < \tau_{2k-4} < 1$, we have
\begin{equation}
\sup_{\| f^{(2k)} \|_{\infty} \le 1} | f(t) - [H_k f] (t) |
\le | S^* (t) - [H_k S^*] (t) | ,
\label{PerfectSplineUpperBound}
\end{equation}
for all $t \in [0,1]$
where $S^*$ is a perfect spline of degree $2k$ with knots $\tau_1, \ldots , \tau_{2k-4}$
which satisfies
\begin{equation}
S^* (t) = (-1)^j, \qquad t \in [\tau_j, \tau_{j+1} )
\label{PerfectSplineCondition}
\end{equation}
for $j= 0, \ldots , 2k-4$ (with the usual convention that $\tau_0 =1$ and $\tau_{2k-3} =1$).
Recall here that a perfect spline $P$ of degree $2k$ with knots
$\tau_1 , \ldots , \tau_{2k-4}$ and satisfying the condition (\ref{PerfectSplineCondition}) is
of the form
$$
P(t) = \sum_{i=0}^{2k-1} \alpha_i t^i + \frac{1}{(2k)!} \left ( t^{2k}
+ 2 \sum_{i=1}^{2k-4} (-1)^i (t-\tau_j)_+^{2k} \right ) \, ,
$$
where
$\alpha_i \in \RR$ for $i=0, \ldots , 2k-1$
(see e.g. [Bojanov, Hakopian, and Sahakian 1993] for a more general form).
The inequality in (\ref{PerfectSplineUpperBound}) says that for any fixed set of
knots the perfect spline gives the largest interpolation error, pointwise, over the
set of $2k-$times continuously differentiable functions such that $\| f^{(2k)} \|_{\infty} \le 1$.
Note that since polynomials of degree $2k-1$ are exactly recovered by the interpolation
operator $H_k$, the perfect spline can be taken to be equal to
$$
S^* (t) = \frac{1}{(2k)!} \left ( t^{2k}
+ 2 \sum_{i=1}^{2k-4} (-1)^i (t-\tau_j)_+^{2k} \right )\, .
$$
Here, we omit writing explicitly the dependence of $S^*$ on the knots, but it should
nevertheless be kept in mind. The bound in (\ref{PerfectSplineUpperBound}) provides
a very useful way to verify computationally the conjecture for the class
$C^{(2k)} [0,1]$ since for any $f \in C^{(2k)} [0,1]$ it follows that
$$
\sup_{0 < \tau_1 < \cdots < \tau_{2k-4}<1}
\| f - [H_k f] \|
\le \| f^{(2k)} \|_{\infty} \sup_{0 < \tau_1 < \cdots < \tau_{2k-4}<1}
\| S^* - [H_k S^*] \|_{\infty} \,,
$$
and hence our conjecture holds for the class $C^{(2k)} [0,1]$ if the right side in
the last display stays bounded.

The bound in (\ref{PerfectSplineUpperBound}) can be proved using for
example the same arguments of the proof of Lemma 6.16 in [Bojanov, Hakopian, and Sahakian 1993].
Indeed, let $ f \in C^{(2k)}[0,1]$ with $\Vert f^{(2k)} \Vert_\infty \leq 1$ and suppose that there exists
$t_0 \in (0,1)$ such that
\begin{eqnarray}\label{Assump}
\vert f(t_0) - [H_k f](t_0) \vert > \vert S^*(t_0) - [H_k S^*](t_0) \vert.
\end{eqnarray}
Note that $t_0 \ne \tau_j, j=0, \cdots, 2k-3$
since $H_k f$ interpolates $f$ at the $\tau_j$'s.
Let
\begin{eqnarray*}
\alpha = \frac{f(t_0) - [H_k f](t_0)}{S^*(t_0) - [H_k S^*](t_0)}
\end{eqnarray*}
and consider the function
\begin{eqnarray*}
h(t) = f(t) - [H_k f](t) - \alpha \left(S^*(t) - [H_k S^*](t)\right), \ \ \ t \in [0,1].
\end{eqnarray*}
By assumption (\ref{Assump}) it follows that $\vert \alpha \vert > 1$.
On the other hand,
the function $h$ admits at least $4k-4 + 1= 4k-3$ zeros: $\tau_j, j=0, \cdots, 2k-2$
(counting multiplicities) and $t_0$. It follows by Rolle's theorem that $h^{(2k-2)}$
has at least $4k-3-(2k-2)= 2k-1$ distinct zeros. Now, note that $h^{(2k-2)}$ is a continuous function
on $[0,1]$ whose first derivative $h^{(2k-1)}$ might only jump at the internal knots.
This implies that $h^{(2k-1)}$ will change its sign at least $2k-2$ times. However,
the latter is impossible. Indeed, since $\vert \alpha \vert > 1$, $\Vert f^{(2k)} \Vert_\infty \leq 1$
and $ \vert (S^*)^{(2k)}(t) \vert = 1, t \in (\tau_j, \tau_{j+1})$ for $j=0, \cdots, 2k-4$, it follows that
\begin{eqnarray*}
\textrm{sign}\ h^{(2k)}(t) = - \textrm{sign}\left(\alpha (S^*)^{(2k)}(t)\right),
\ \ \ \ \textrm{for} \ t \in (\tau_j, \tau_{j+1})
\end{eqnarray*}
and therefore the function $h^{(2k)}$ will have the same number of sign changes as
$(S^*)^{(2k)}$, namely $2k-4$, and the sign changes occur when $t$ takes values in the set of
internal knots $\tau_j$, $j=1, \ldots , 2k-4$.
This in turn implies that $h^{(2k-1)}$ has at most $2k-3$ sign changes. The contradiction completes the proof.

Our numerical results based on $10^4$ Monte Carlo replications
point to boundedness of the interpolation error for the perfect spline
and are reported in Table \ref{PerfectSplineErrors} for $k=3,\cdots, 10$ where the column labelled ``factor''
multiplies all of the preceding columns.
Note the particular pattern of decay of the maximal error as $k$ increases.
\begin{center}
\begin{table}[12pt]
\caption{Statistical summary of the sup-norm of the interpolation error, $\mathcal{E}_k(S^*)$, for $k=3,\cdots,10$ based on $10^4$ simulations. Table produced using the Mathematica program \lq PS-MC-post\rq.}
\medskip
\begin{tabular}{|c|cccccccc|}
\hline
$k$ & Mean & Std. dev. & Median & $95\%$-q & $99\%$-q & max & factor & max $\times (2k)!$\\
\hline\hline
3 & 1.71 & 1.10 & 2.15 & 2.77 & 2.78 & 2.78 & $\times 10^{-3}$ & 2.001\\ \hline
4 & 2.77 & 1.86 & 3.28 & 4.94 & 4.96 & 4.96 & $\times 10^{-5}$ & 1.999 \\ \hline
5 & 2.84 & 2.04 & 3.18 & 5.46 & 5.50 & 5.51 & $\times 10^{-7}$ & 1.999\\ \hline
6 & 2.08 & 1.53 & 2.33 & 4.12 & 4.17 & 4.18 & $\times 10^{-9}$ & 2.002 \\ \hline
7 & 1.09 & 8.37 & 1.19 & 2.25 & 2.29 & 2.29 & $\times 10^{-11}$ & 1.996 \\ \hline
8 & 4.36 & 3.51 & 5.03 & 8.99 & 9.52 & 9.54 & $\times 10^{-14}$ & 1.996\\ \hline
9 & 1.47 & 1.08 & 1.43 & 3.00 & 3.08 & 3.10 & $\times 10^{-16}$ & 1.987 \\ \hline
10 & 3.86 & 2.94 & 4.34 & 8.00 & 8.19 & 8.21 & $\times 10^{-19}$ & 1.999 \\ \hline
\end{tabular}
\label{PerfectSplineErrors}
\end{table}
\end{center}

\bigskip

The last column corresponds to the maximal error for interpolating
the perfect spline $\times (2k)!$. This gives an upper bound for
the error for interpolating $f_{0}(t)= t^{2k}$; i.e., an upper bound
for the supremum norm of the monospline $f_0 - H_k f_0$.
The obtained products seem to remain very close to 2, which
is rather an interesting outcome in its own right. We should
mention here that simulations with a larger number of replications
($2 \times 10^4$ and $5 \times 10^4$) yielded similar statistics
to those reported in Table \ref{PerfectSplineErrors}.\\

\begin{center}
\begin{table}[12pt]
\caption{Statistical summary of the sup-norm of the interpolation
error given by the Hermite (first row) and complete (second row)
spline interpolants for interpolating $f_0(t) = t^{2k}$ and $k=3,\cdots,10$
based on $10^4$ simulations. Table produced using the Mathematica programs \lq MN-MC-post\rq and \lq MS-MC-ConfigIsol-post\rq.}
\medskip
\begin{tabular}{|c|cccccc|}
\hline
$k$ & Mean & Median & Std. dev. & $95\%$-q & $99\%$-q & max \\
\hline
\hline
3 & $2.52 \times 10^{-3}$ & $1.17\times 10^{-3}$ & $3.24 \times 10^{-3}$ & $9.68 \times 10^{-3}$ & $1.47 \times 10^{-2} $ & $2.14 \times 10^{-2}$ \\
3 & $2.74 \times 10^{-3}$ & $1.62\times 10^{-3}$ & $2.85 \times 10^{-3}$ & $8.93 \times 10^{-3}$ & $1.22 \times 10^{-2} $ & $1.53 \times 10^{-2}$ \\
\hline
4 & $1.36 \times 10^{-4}$ & $4.18\times 10^{-5}$ & $3.12 \times 10^{-4}$ & $5.44 \times 10^{-4}$ & $1.5 \times 10^{-3} $ & $6.89 \times 10^{-3}$ \\
4 & $1.27 \times 10^{-4}$ & $4.14\times 10^{-5}$ & $2.39 \times 10^{-4}$ & $5.63 \times 10^{-4}$ & $1.22 \times 10^{-3} $ & $3.24 \times 10^{-3}$ \\
\hline
5 & $1.67 \times 10^{-5}$ & $3.2 \times 10^{-6}$ & $10^{-4}$ & $6.67 \times 10^{-5}$ & $2 \times 10^{-4}$ & $8.5 \times 10^{-3}$ \\
5 & $5.74 \times 10^{-6}$ & $1.04 \times 10^{-6}$ & $1.64 \times 10^{-5}$ & $2.66 \times 10^{-5}$ & $7.31 \times 10^{-5}$ & $3.97 \times 10^{-4}$ \\
\hline
6 & $2.26 \times 10^{-4}$ & $4.13 \times 10^{-7}$ & $1.1 \times 10^{-2}$ & $3.31 \times 10^{-5}$ & $2.45 \times 10^{-4}$ & 1.02 \\
6 & $2.6 \times 10^{-7}$ & $2.51 \times 10^{-8}$ & $1.09 \times 10^{-6}$ & $1.07 \times 10^{-6}$ & $4.55 \times 10^{-6}$ & $3.77 \times 10^{-5}$ \\
\hline
7 & $4.01$ & $8.63 \times 10^{-8}$ & $4.01 \times 10^2$ & $2.93 \times 10^{-5}$ & $7.97 \times 10^{-4}$ & $4.01 \times 10^{5}$ \\
7 & $1.16 \times 10^{-8}$ & $5.52 \times 10^{-10}$ & $8.35 \times 10^{-8}$ & $3.91 \times 10^{-8}$ & $2.02 \times 10^{-7}$ & $5.21 \times 10^{-6}$ \\
\hline
8 & $1.83 \times 10^{6}$ & $3.13 \times 10^{-8}$ & $1.83 \times 10^{8}$ & $4.35 \times 10^{-5}$ & $4.43 \times 10^{-3}$ & $1.83 \times 10^{10}$ \\
8 & $5.93 \times 10^{-10}$ & $1.41 \times 10^{-11}$ & $6.94 \times 10^{-9}$ & $1.50 \times 10^{-9}$ & $1.15 \times 10^{-8}$ & $5.65 \times 10^{-7}$ \\
\hline
9 & $1.15 \times 10^{5}$ & $1.7 \times 10^{-8}$ & $1.14 \times 10^{6}$ & $1.5 \times 10^{-4}$ & $3.8 \times 10^{-2}$ & $1.14 \times 10^{8}$ \\
9 & $2.07 \times 10^{-11}$ & $3.03 \times 10^{-13}$ & $2.53 \times 10^{-10}$ & $4.74 \times 10^{-11}$ & $2.81 \times 10^{-10}$ & $1.65 \times 10^{-8}$ \\
\hline
10 & $8.43 \times 10^{10}$ & $1.19 \times 10^{-8}$ & $8.42 \times 10^{12}$ & $4.05 \times 10^{-4}$ & $2.12 \times 10^{-2}$ & $8.42 \times 10^{14}$ \\
10 & $9.65 \times 10^{-13}$ & $6.51 \times 10^{-15}$ & $2.26 \times 10^{-11}$ & $1.41 \times 10^{-12}$ & $1.46 \times 10^{-11}$ & $2.11 \times 10^{-9}$ \\ \hline
\end{tabular}
\label{Monospline}
\end{table}
\end{center}

A natural question to ask is whether the obtained bounds are actually
achieved by the corresponding monosplines. Based on $10^4$ simulations,
a statistical summary of
$\Vert f_0 - H_k f_0 \Vert_\infty = \Vert \mathcal{E}_k(f_0) \Vert $
and that of the complete interpolation error are obtained for $k=3, \cdots, 10$
and reported in Table \ref{Monospline}. The results show two different aspects.
On one hand, we see that the maximal Hermite interpolation error is much
smaller than 2 for the values $k=3,4,5$, and that the errors given by the
Hermite and complete spline interpolants are quite comparable.
On the other hand, when $k \geq 7$, the obtained maximal errors
given by the Hermite spline are very large. We believe that these large
values
result from numerical instabilities as they are not at all in
agreement with the rather stable upper bounds given by the perfect splines.
In order to investigate more this latter aspect, the 10 largest interpolation
errors have been isolated but more importantly the corresponding
configurations of the knots (for more numerical details,
see the Mathematica program \lq MN-MC-ConfigIsol-post\rq).
In Table \ref{KnotConfig} we report the largest and second largest error
that have occurred over $10^4$ independent replications and the
corresponding knot configurations, for $k=7,\cdots, 10$. In the table,
we underline the knots whose in-between distances are of the order of $10^{-3}$.
The closeness of two knots or more is expected to result in severe ill-conditioning
of the linear system to be solved. For these \lq\lq bad\rq\rq \ configurations,
stability of the complete problem might be partially explained by the fact that
only the function is interpolated. As we interpolate also the first derivatives
at the knots, having them very close to each other makes the linear
system in our problem much closer to being singular.

\begin{center}
\begin{table}[12pt]
\caption{Knot configurations for the largest and second largest maximal
Hermite interpolation errors ($e_1$ and $e_2$) for interpolating $f_0$
and $k=7,8,9,10$ based on $10^4$ simulations. Table produced using the Mathematica program
\lq MS-MC-ConfigIsol-post\rq.}
\medskip
\begin{tabular}{|c|c|c||c|c|}
\hline
$k$ & $e_1$ & Config. & $e_2$ & Config. \\
\hline
\hline
& & 0.022 \ 0.065 \ 0.095 & & \underline{0.3433} \ \underline{0.3439} \ \underline{0.3444} \\
& & \underline{0.269} \ \underline{0.272} \ 0.377 & & 0.4709 \ 0.5058 \ 0.5327 \\
7 & $4.01 \times 10^5$ & 0.582 \ \underline{0.678} \ \underline{0.685} \ \underline{0.686} & 2.05 & 0.6057 \ \underline{0.9460} \ \underline{0.9472} \ 0.999 \\
\hline
& & 0.0810 \ 0.1265 \ 0.1360 & & 0.0178 \ 0.0709 \ 0.0960 \\
& & 0.1410 \ 0.1573 \ 0.1680 & & 0.105 \ 0.153 \ 0.279 \\
& & 0.3770 \ 0.3820 \ 0.6975 & & 0.308 \ 0.366 \ 0.380 \\
8 & $1.83 \times 10^{10}$ & \underline{0.7873} \ \underline{0.7873} \ \underline{0.7879} & $7.09 \times 10^6 $ & \underline{0.7477} \ \underline{0.7478} \ 0.7505 \\
\hline
& & \underline{0.2901} \ \underline{0.2905} \ \underline{0.2933} \ 0.4046 & & 0.023 \ 0.064 \ 0.128 \ 0.160 \\
& & 0.442 \ 0.510 \ 0.530 & & 0.198 \ 0.275 \ 0.399 \\
& & \underline{0.732} \ \underline{0.747} \ 0.891 & & 0.551 \ 0.593 \ 0.624 \\
9 & $1.14 \times 10^{8}$ & 0.894 \ \underline{0.903} \ \underline{0.943} \ 0.973 & $7.71 \times 10^5 $ & 0.652 \ \underline{0.683} \ \underline{0.686} \ \underline{0.688} \\
\hline
& & \underline{0.241} \ \underline{0.241} \ 0.292 \ 0.322 & & \underline{0.1721} \ \underline{0.1726} \ \underline{0.1739} \ 0.2551 \\
& & 0.340 \ 0.362 \ 0.375 \ 0.436 & & 0.382 \ 0.601 \ 0.613 \ 0.701 \\
& & 0.467 \ 0.586 \ 0.704 \ 0.761 & & \underline{0.743} \ \underline{0.745} \ 0.763 \ 0.783 \\
10 & $ 8.42 \times 10^{14} $& 0.789 \ \underline{0.8890} \ \underline{0.8892} \ \underline{0.8894} & $ 6.17 \times 10^{11} $ & 0.802 \ 0.826 \ 0.953 \ 0.965 \\ \hline
\end{tabular}
\label{KnotConfig}
\end{table}
\end{center}
Although the following fact is rather
peripheral to the main problem treated
in this paper, we would like to record it here. The result was
brought to our attention by C. de Boor.
When the knots of the spline interpolant are equally spaced; i.e.,
$\tau_{j+1} - \tau_j = \frac{1}{2k-3}, j =0, \cdots, 2k-4$, the interpolation
error $\mathcal{E}_k(f_0)$ is the periodic monospline of degree $2k$
of period $1/(2k-3)$, that we denote here by $M_{2k}$. One can prove
that this monospline is of one sign on $[0,1]$ as it has $4k-4$ (double)
zeros, and that this number is maximal [Corollary 1, page 422 of Micchelli 1972;
Theorem 7.1 in Bojanov, Hakopian and Sahakian 1993].
Furthermore, it can be shown that, on $[0,1/(2k-3)]$, we have
\begin{eqnarray*}
M_{2k}(t) = \frac{1}{(2k-3)^{2k}} \big (\mathcal{B}_{2k}((2k-3)t) - B_{2k} \big),
\end{eqnarray*}
where $\mathcal{B}_{2k}$ is the Bernoulli polynomial of degree $2k$,
and $B_{2k}=\mathcal{B}_{2k}(0)$ is the corresponding Bernoulli number.
In this case, the $L_\infty$ norm of the interpolation error is given by
\begin{eqnarray*}
\Vert M_{2k} \Vert_\infty &=& \frac{1}{(2k-3)^{2k}} \big \vert \mathcal{B}_{2k}(1/2) -B_{2k} \big \vert \\
&=& \frac{2}{(2k-3)^{2k}} \big \vert (1- 2^{-2k}) \times B_{2k} \big \vert,
\end{eqnarray*}
using the identity $\mathcal{B}_{2k}(1/2) = - (1-2^{1-2k}) B_{2k}$.
The order of the maximal error decreases quickly with $k$, and its
values for $3 \leq k \leq 6$ are reported in Table \ref{MaxMonoS}.
We conjecture that, for a fixed $k$, this maximal error is the
smallest over all configurations of the knots.

\begin{center}
\begin{table}[12pt]
\caption{Some values of the $L_\infty$-norm of the interpolation error
for interpolating $t^{2k}$ on $[0,1]$ when the knots are equispaced. Table produced using the Mathematica program \lq MN-SinglePrint-post\rq.}
\medskip
\begin{tabular}{|c|c|c|}
\hline
$k$ & Bernoulli number $B_{2k}$ & Maximal error \\
\hline
\hline
3 & $1 /42$ & $1 /15552 \approx 6.3 \times 10^{-5}$ \\
\hline
4 & $ -1/30$ & $17/100000000 = 1.7 \times 10^{-7}$ \\
\hline
5 & $ 5/66 $ & $155/289254654976 \approx 5.35\times 10^{-10}$\\
\hline
6 & $-691/2730$ & $691/385610460475392 \approx 1.71 \times 10^{-12}$\\ \hline

\end{tabular}
\label{MaxMonoS}
\end{table}
\end{center}

\section{Conclusions and open questions}

As mentioned in section 2, our Hermite interpolation
arises naturally in the
study of the stochastic behavior
of the
Least Squares and Maximum Likelihood estimators of a $k$-monotone density.
Considering interpolation of the hinge functions $f_u, u \in (0,1)$ might
appear somewhat too ambitious as boundedness of the error
independently of the locations of the knots in this problem is only
a sufficient condition for the result to hold for smoother classes of functions.
But as a matter of fact, hinge functions constitute a basis of the $(k-1)$-fold
integral of empirical distribution functions which are involved in the original
estimation problem and play a major a role in understanding the asymptotic
properties of the estimators (for more technical details see
[Balabdaoui and Wellner 2004d, page 16, Lemma 2.4]).
Our numerical investigations point partially to the validity of the
conjecture for small values of $k$ but leave us with more unanswered
questions for larger values.

\begin{figure}[H]
\begin{center}
\includegraphics[bb=63 1 459 236, clip=true, width=11cm]{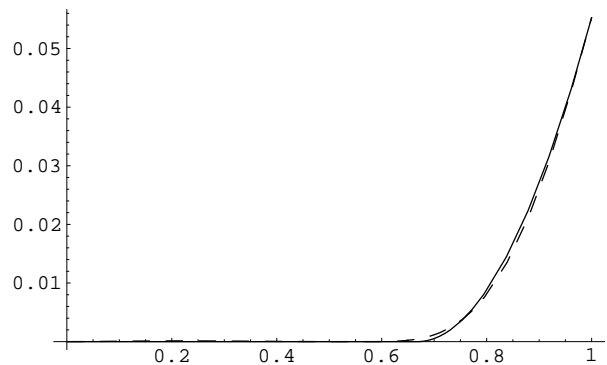}
\end{center}
\caption{Plot of the hinge function $f_u(t)= (t-u)^{k-1}_+ /(k-1)!$ with $u= 0.66$ and
$k=3$ (solid line) and its Hermite spline interpolant (dashed line). The interior knots are 0.54, 0.55. Plot produced using the Mathematica program \lq EB-SinglePrint-hinge-post\rq.}
\label{Hinge-Hinterp-K3}
\end{figure}

\begin{figure}[H]
\begin{center}
\includegraphics[bb=69 17 459 236, clip=true, width=11cm]{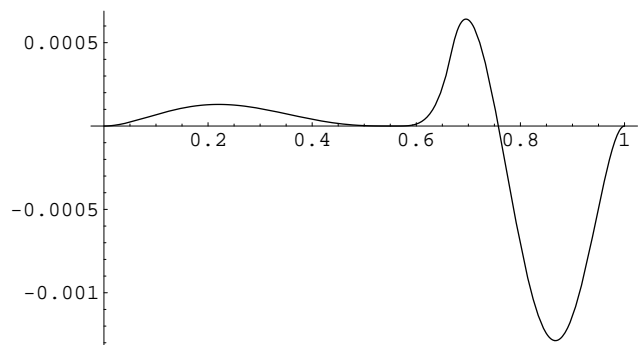}
\end{center}
\caption{Plot of the Hermite interpolation error for interpolating the hinge function
$f_u(t)= (t-u)^{k-1}_+ /(k-1)!$ with $u= 0.66$ and $k=3$. The interior knots are 0.54, 0.55.
The maximal error is $1.28 \times 10^{-3}$. Plot produced using the 
Mathematica program \lq EB-SinglePrint-hinge-post\rq.}
\label{Hinge-Herror-K3}
\end{figure}

On the other hand, our simulations based on the perfect splines theory
suggest strongly that the conjecture is true when the interpolated function
belongs to the $C^{(2k)}[0,1]$, another class that is also of great relevance
for us in connection with the same estimation problems (see
[Balabdaoui and Wellner 2004d, page 16, Lemma 2.4]). In the particular
case of $f(t) = f_0(t) = t^{2k}$, it follows from the maximal error for
interpolating the perfect spline that there exists some constant $d_k > 0$
very close to 2 such that
\begin{eqnarray*}
\sup_{0 < \tau_1 < \cdots < \tau_{2k-4} < 1}
\Vert f_0 - H_k f_0 \Vert_\infty \leq d_k \Vert f^{(2k)} \Vert_\infty
\end{eqnarray*}
The importance of the latter result is two-fold: beside that it gives
a uniform upper bound on the monospline (it turns out that this
bound can be greatly improved for small values of $k$, see
Table \ref{Monospline}), it serves as a way of checking the
numerical stability of the solution given by interpolating
\lq\lq directly\rq\rq \ $f_0$. The extremely large values obtained
for the $\Vert \mathcal{E}_k(f_0) \Vert_\infty$ for $k=7, \cdots, 10$
are believed to be due to severe ill-conditioning of the problem as
some of the knots are very close to each other (see Table \ref{KnotConfig}).
This might explain also the large errors obtained for interpolating the
hinge functions $f_u, u \in (0,1)$ for large values if $k$.
\begin{figure}[H]
\begin{center}
\includegraphics[bb=67 1 459 236, clip=true, width=11cm]{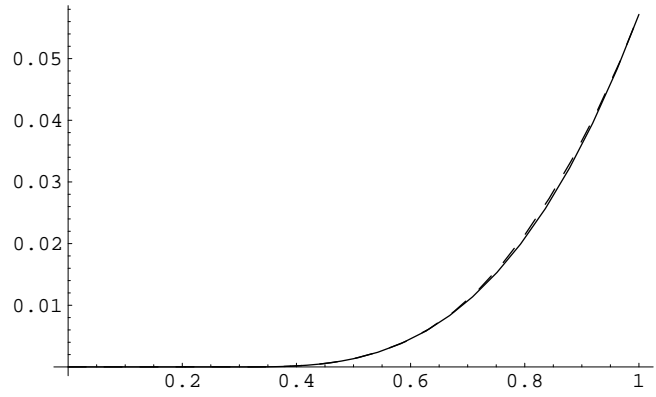}
\end{center}
\caption{Plot of the hinge function $f_u(t)= (t-u)^{k-1}_+ /(k-1)!$ with $u= 0.30$ and
$k=4$ (solid line) and its Hermite spline interpolant (dashed line).
The interior knots are 0.11, 0.33, 0.49, 0.50. Plot produced using the Mathematica program \lq EB-SinglePrint-hinge-post\rq.}
\label{Hinge-Hinterp-K4}
\end{figure}
\begin{figure}[H]
\begin{center}
\includegraphics[bb=67 1 459 236, clip=true, width=11cm]{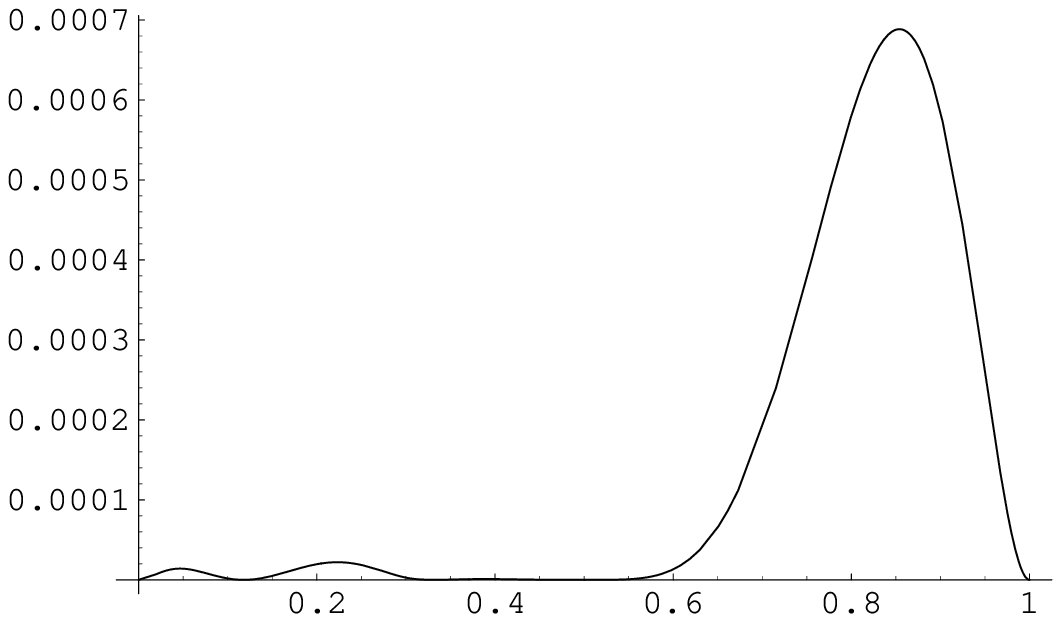}
\end{center}
\caption{Plot of the Hermite interpolation error for interpolating the
hinge function $f_u(t)= (t-u)^{k-1}_+ /(k-1)!$ with $u= 0.30$ and $k=4$.
The interior knots are $0.11, 0.33, 0.49, 0.50$. The maximal error is
$6.8 \times 10^{-4}$. Plot produced using the Mathematica program \lq EB-SinglePrint-hinge-post\rq.}
\label{Hinge-Herror-K4}
\end{figure}

\begin{figure}[H]
\begin{center}
\includegraphics[bb=69 17 459 236, clip=true, width=11cm]{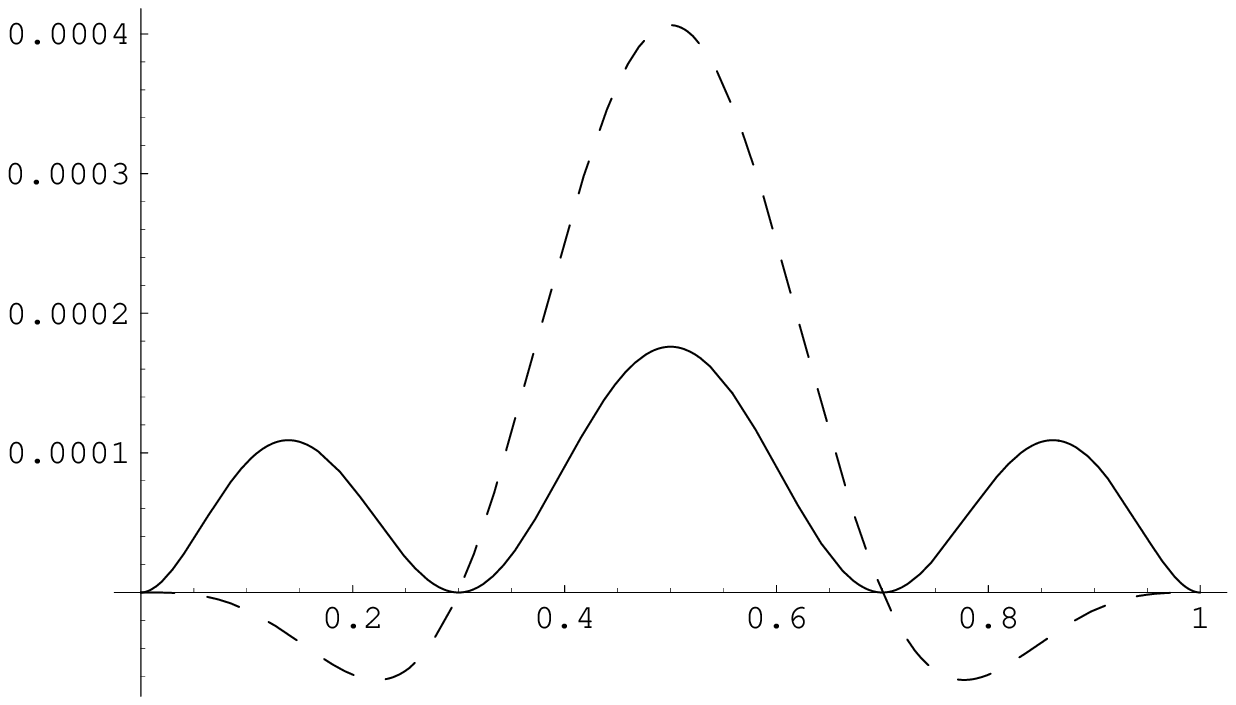}
\end{center}
\caption{Plot of the Hermite (solid line) and complete (dashed line)
interpolation error for $f_0(t)=t^{2k}$ on $[0,1]$ for $k=3$.
The maximal Hermite and complete interpolation errors are
$1.76 \times 10^{-4}$ and $4 \times 10^{-4}$ respectively.
The interior knots are $0.30$, $0.70$. Plot produced using the Mathematica program \lq MN-SinglePrint-post\rq.}
\label{MN-HC-K3}
\end{figure}

\begin{figure}[H]
\begin{center}
\includegraphics[bb=69 17 459 236, clip=true, width=11cm]{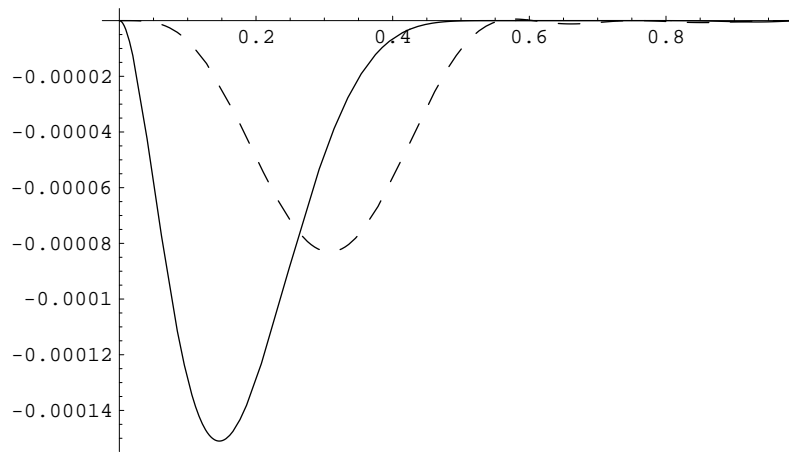}
\end{center}
\caption{Plot of the Hermite (solid line) and complete (dashed line)
interpolation error for $f_0(t)=t^{2k}$ on $[0,1]$ for $k=4$.
The maximal Hermite and complete interpolation errors are
$1.5 \times 10^{-4}$ and $8.31 \times 10^{-5}$ respectively.
The interior knots are 0.55, 0.60, 0.74, 0.76. Plot produced using the Mathematica program \lq MN-SinglePrint-post\rq.}
\label{MN-HC-K4}
\end{figure}

A part of this experimental work was not presented here and deals
mainly with further computations aiming at checking whether there
is any effect of using the B-spline basis as opposed to the canonical basis.
The outcomes of many simulations that were performed with the
two different bases and taking for every replication the same random
knots were very comparable. On the other hand, we have done more
comparisons between our Hermite interpolant and the complete
one to get a better understanding about their relative goodness
of approximation. One might think that our spline interpolant
problem should behave better than the complete spline in the
case where the number of (interior) knots $m-1 =2k-4$: the
reasoning would be based on the fact that our spline interpolant
matches not only the function to be approximated at the knots
(including the boundary points)
but it also matches its tangent, and hence tries to ``adapt'' \
more to its shape. In all the examples that we have taken
including $f_0(t) = t^{2k}$ disprove this idea.
But the question remains open: Are there functions that are better
approximated (in the $L_\infty$ sense) by our Hermite spline than
by the complete spline$?$ This question prompts another one in a
different direction: We know that the perfect spline gives in our
problem the largest interpolation error over the class $C^{(2k)}[0,1]$
and even over Sobolev space
$W_{\infty}^{2k}$, but what kind of functions give the worst interpolation
error in the Sobolev space $W_{\infty}^{k-1}?$ In the complete
interpolation problem, [Shadrin 2001] makes use of the properties
of the null spline to bound the orthogonal projector of $C[0,1]$ equipped
with the $L_\infty$-norm on the space of splines of degree $k-1$.
The argument in our problem might need to involve such constructions
as well. Finally, it seems natural to wonder whether, for a given
sufficiently smooth function $f$, the spline $H_kf$ is the solution
of some minimization problem: We wonder whether the $L_2$-norm
(or in general the $L_p$-norm, $p \geq 2$) of $(H_kf)^{(k)}$
(or in general of $(H_kf)^{(j)}$ for some $j \leq k$) minimizes some
functional, yet to be defined. The problem seems to be hard to solve,
and it is not clear yet in which direction one should try to find such a
criterion, although the interpolation problem was first originated by
a Least Squares problem.\\

Our aim in this paper is to call attention
to this new interpolation problem, which plays the role of a bridge between
deterministic spline theory and statistical estimation for $k-$monotone
functions.

\bigskip

\par\noindent
{\bf Acknowledgments:}
We would like to thank Carl de Boor, Nira Dyn,
Grace Wahba, and Alexi Yu. Shadrin for sharing with us their thoughts and
suggestions about this interpolation problem, and Tilmann Gneiting
for several helpful suggestions. Discussions with Moulinath Banerjee and
Marloes Maathuis during a RiP at Oberwolfach led to serveral improvements
and additions.

\end{document}